\documentclass{amsart}
\usepackage{amscd,amssymb}
\usepackage[all]{xy}

\newcommand{\gr}{{\operatorname{gr}\nolimits}}

\newcommand{\Hom}{\operatorname{Hom}\nolimits}

\renewcommand{\Im}{\operatorname{Im}\nolimits}

\newcommand{\rrad}{\mathfrak{r}}

\newcommand{\Ext}{\operatorname{Ext}\nolimits}

\newcommand{\op}{{\operatorname{op}\nolimits}}

\newcommand{\id}{{\operatorname{id}\nolimits}}

\newcommand{\fraka}{\mathfrak{a}}

\newcommand{\frako}{\mathfrak{o}}
\newcommand{\frakt}{\mathfrak{t}}

\renewcommand{\L}{\Lambda}
\newcommand{\Z}{{\mathbb Z}}

\newcommand{\extto}{\xrightarrow}

\newcommand{\HH}{\operatorname{HH}\nolimits}

\newtheorem{lem}{Lemma}[section]
\newtheorem{prop}[lem]{Proposition}

\newtheorem{thm}[lem]{Theorem}
\newtheorem*{Theorem}{Theorem}
\theoremstyle{definition}

\newtheorem{example}[lem]{Example}

\begin{document}

\title{Multiplicative structures for Koszul algebras}

\author[Buchweitz]{Ragnar-Olaf Buchweitz}
\address{Ragnar-Olaf Buchweitz\\ Department of Mathematics\\
University of Toronto\\ To\-ronto, ON Canada M5S 3G3\\ Canada}
\email{ragnar@math.toronto.edu}
\author[Green]{Edward L. Green$^\dagger$}
\thanks{$^\dagger$Partially supported by a grant from the NSA} 
\address{Edward L. Green\\ Department of Mathematics\\
Virginia Tech\\ Blacksburg, VA 24061\\ USA}
\email{green@math.vt.edu}
\author[Snashall]{Nicole Snashall}
\address{Nicole Snashall\\ 
Department of Mathematics\\
University of Leicester\\
University Road\\
Leicester, LE1 7RH\\
England}
\email{N.Snashall@mcs.le.ac.uk}
\author[Solberg]{\O yvind Solberg$^\ddagger$}
\thanks{$^\ddagger$Partially supported by the
  Research Council of Norway}
\address{\O yvind Solberg\\Institutt for matematiske fag\\
NTNU\\ N--7491 Trondheim\\ Norway}
\email{oyvinso@math.ntnu.no}

\date{\today}

\begin{abstract}
Let $\L=kQ/I$ be a Koszul algebra over a field $k$, where $Q$ is a
finite quiver. An algorithmic method for finding a minimal projective
resolution $\mathbb{F}$ of the graded simple modules over $\L$ is
given in \cite{GS}. This resolution is shown to have a
``comultiplicative'' structure in \cite{GHMS}, and this is used to
find a minimal projective resolution $\mathbb{P}$ of $\L$ over the
enveloping algebra $\L^e$. Using these results we show that the
multiplication in the Hochschild cohomology ring of $\L$ relative to
the resolution $\mathbb{P}$ is given as a cup product and also provide
a description of this product. This comultiplicative structure also
yields the structure constants of the Koszul dual of $\L$ with respect
to a canonical basis over $k$ associated to the resolution
$\mathbb{F}$. The natural map from the Hochschild cohomology to the
Koszul dual of $\L$ is shown to be surjective onto the graded centre
of the Koszul dual.
\end{abstract}

\maketitle

\section*{Introduction}
The aim of this paper is to show that the knowledge of a minimal
projective resolution of the graded simple modules for a Koszul
algebra over a field $k$ not only provides a $k$-basis for the Koszul
dual with structure constants, but also gives a closed formula for the
multiplication in the Hochschild cohomology ring. In a forthcoming
paper \cite{GS} an algorithm will be presented for computing a minimal
projective resolution of linear modules for a Koszul algebra, and thus
making Koszul algebras more computationally accessible.

A prime example of a Koszul algebra is a polynomial ring in a finite
number of commuting variables. Its Koszul dual is the exterior
algebra, which is a finite dimensional selfinjective algebra. Using
techniques from representation theory of finite dimensional algebras
and Koszul duality, Mart\'inez-Villa and Zacharia have obtained new
results on the structure of locally free sheaves on projective space
\cite{MZ}. Polynomial rings are Artin-Schelter regular Koszul algebras
and their Koszul duals are finite dimensional selfinjective
algebras. More generally, this connection between representation
theory of finite dimensional algebras and Artin-Schelter regular
Koszul algebras allows one to apply the theory of support varieties
for finite dimensional selfinjective algebras defined via the
Hochschild cohomology ring (see \cite{EHSST,SS}). To apply the theory
of support varieties some finiteness conditions must be
satisfied. These conditions are related to the natural map from the
Hochschild cohomology ring to the Koszul dual. The image of this map
is characterized in Section \ref{section:4} of this paper.

We now introduce notation we use throughout the paper and describe the
results of the paper in more detail.  Let $Q$ be a finite quiver and
$k$ a field. We denote the ideal generated by the arrows in $Q$ in the
path algebra $kQ$ by $J$.  Let $\L=kQ/I$ be a Koszul algebra, where
$I$ is an ideal contained in $J^2$. Set $\rrad=J/I$.  In \cite{GSZ} a
method for constructing a minimal projective resolution of $\L/\rrad$
as a right $\L$-module was presented. If $R=kQ$, then this minimal
right resolution of the degree zero part $\L_0$ of $\L$ can be given
in terms of a filtration of right ideals
\[\cdots\subseteq \amalg_{i=0}^{t_n} f^n_iR \subseteq
\amalg_{i=0}^{t_{n-1}}  
f^{n-1}_iR\subseteq \cdots \subseteq \amalg_{i=0}^{t_1}
f^1_iR\subseteq \amalg_{i=0}^{t_0}f^0_iR=R\] in $R$.  In \cite{GHMS}
it is shown that there is a comultiplicative structure associated to the
$\{f^n_i\}$; in particular there exist $c_{pq}(n,i,r)$ in $k$ such
that
\begin{equation}
f^n_i  = \sum_{p=0}^{t_r}\sum_{q=0}^{t_{n-r}}
c_{pq}(n,i,r)f^r_pf^{n-r}_q\label{eq:comult}
\end{equation}
for all $n\geq 1$, all $i$ in $\{0,1,\ldots,t_n\}$ and all $r$ in
$\{0,1,\ldots,n\}$. 

The paper begins with recalling a reduced bar resolution introduced in
\cite{C2}. As in the standard setup \cite{Ge}, there is an associated
cup product, which we show is the same as the Yoneda product in the
Hochschild cohomology ring. This reduced bar resolution together with
the comultiplicative structure is shown to provide a description of
the multiplication in the Hochs\-child cohomology ring $\HH^*(\L)$ of
$\L$ in Section \ref{section:2}.  Recall that the enveloping algebra
$\L^e$ of $\L$ is given by $\L^\op\otimes_k \L$. In our description of
the multiplicative structure of $\HH^*(\L)$, we use the minimal
projective resolution of $\L$ as a right $\L^e$-module presented in
\cite{GHMS}.

Recall that an element $x$ in $kQ$ is called \emph{uniform} if $x$ is
non-zero and there exist vertices $u$ and $v$ in $Q$ such that
$x=uxv$. If $x$ is a uniform element with $x=uxv$ with $u$ and $v$
vertices, then we write $\frako(x)=u$ and $\frakt(x)=v$. It was shown
in \cite{GSZ} that the $\{f^n_i\}$ can be chosen to be uniform
elements.  The projective modules occurring in a minimal projective
resolution of $\L$ over $\L^e$ are given by $P^n=\amalg_{i=0}^{t_n}
\L\frako(f^n_i)\otimes_k \frakt(f^n_i)\L$ for $n\geq 0$ (see
\cite{H}). In Section \ref{section:2} the multiplicative structure of
the Hochschild cohomology ring of a Koszul algebra $\L$ is shown to
have the following description.

\begin{Theorem}
Suppose that $\eta\colon P^n\to \L$ and $\theta\colon P^m\to \L$
represent elements in $\HH^*(\L)$ and are given by
$\eta(\frako(f^n_i)\otimes \frakt(f^n_i))=\lambda_i$ for
$i=0,1,\ldots,t_n$ and $\theta(\frako(f^m_i)\otimes
\frakt(f^m_i))=\lambda_i'$ for $i=0,1,\ldots,t_m$. Then
\[\eta\ast\theta(\frako(f^{n+m}_i)\otimes \frakt(f^{n+m}_i))
= \sum_{p=0}^{t_n}\sum_{q=0}^{t_m}c_{pq}(n+m,i,n)
\lambda_{p}\lambda_{q}'\]
for $i=0,1,\ldots,t_{m+n}$, where $\eta\ast\theta$ represents the
product of $\eta$ and $\theta$ in $\HH^*(\L)$.
\end{Theorem}

The multiplicative structure of the Hochschild cohomology ring is not
known for most algebras, including Koszul algebras.  One class of
Koszul algebras for which the product is known is radical square zero
algebras $\L=kQ/J^2$ with $Q$ not an oriented cycle; then all products
of elements in $\HH^{\geq 1}(\L)$ are zero by \cite{C}. Another class
of Koszul algebras where the ring structure has been determined can be
found in \cite{BGMS}, where, for algebras in this class, $\HH^*(\L)$
is a finitely generated algebra.

Since $\L$ is a Koszul algebra, its Koszul dual,
$E(\L)=\Ext^*_\L(\L/\rrad,\L/\rrad)$, is also a Koszul algebra.  In
Section \ref{section:3} we show that the comultiplicative structure of
the sets $\{ f^n_i\}$ can be used to construct the structure constants
of the Koszul dual of $\L$ for a basis associated to the sets 
$\{f^n_i\}$.  In particular, given the sets $\{f^n_i\}$, we construct a
basis $\{\widehat{f}^n_i\}$ of $E(\L)$ so that the product is given by
\[\widehat{f}^m_i\widehat{f}^n_j=\sum_{l=0}^{t_{m+n}}
c_{ji}(m+n,l,n)\widehat{f}^{m+n}_l.\]
See Theorem \ref{thm:structconst} for details. 

There is a natural map $\varphi_{\L/\rrad}\colon \HH^*(\L)\to E(\L)$
whose image is contained in the graded centre of $E(\L)$. Recall that
the graded centre $Z_\gr(E(\L))$ is the subring of $E(\L)$ generated by
all homogeneous elements $z$ such that $zy=(-1)^{|z||y|}yz$ for each
homogeneous element $y$ in $E(\L)$, where $|x|$ denotes the degree of
a homogeneous element $x$. In Section \ref{section:4} we apply the
results of the previous sections to show that $\varphi_{\L/\rrad}$ has
image $Z_\gr(E(\L))$. 

The paper ends with a number of examples illustrating our results. 

Throughout, we consider right modules unless otherwise explicitly
stated.

\section{A reduced bar resolution}

For this section we fix the following notation and assumptions. Let
$\L$ be an algebra over a field $k$, and let $\fraka$ be an ideal in
$\L$. Set $\L_0=\L/\fraka$ and assume that the natural algebra
homomorphism $\L\to \L_0$ is a split $k$-algebra
homomorphism. Finally, assume that $\L_0$ is a finite product of
copies of $k$. For example, if $\L=kQ/I$ with $I$ contained in $J^2$
and $\fraka=J/I$, then $\L$ satisfies these conditions. 

In this section we recall a reduced bar resolution $(\mathbb{B},d)$ of
$\L$ introduced in \cite{C2}. If $\L_0=\prod_{i=1}^m k$, then we let
$\{e_1,\ldots,e_m\}$ be a complete set of primitive orthogonal central
idempotents of $\L_0$.  Define $(\mathbb{B},d)$ by
$B^n=\L^{\otimes_{\L_0}(n+2)}$, the $(n+2)$-fold tensor product of
$\L$ over $\L_0$, and $d^n\colon B^n\to B^{n-1}$ by
\[d^n(\lambda_0\otimes \cdots\otimes \lambda_{n+1})=
\sum_{i=0}^n(-1)^i(\lambda_0\otimes \cdots\otimes
\lambda_i\lambda_{i+1}\otimes \cdots\otimes \lambda_{n+1}).\] It is
shown in \cite[Lemma 1.1]{C2} that $(\mathbb{B},d)$ is a projective
$\L^e$-resolution of $\L$. When $\L_0=k$, this is the usual bar
resolution. Note that if $\L_0=\prod_{i=1}^m k$ for $m>1$, then $\L_0$
is not central in general, so that $\L$ is not necessarily an algebra
over $\L_0$.

For the bar resolution $(\mathbb{B},d)$ there is a chain map
$\Delta\colon \mathbb{B}\to \mathbb{B}\otimes_\L\mathbb{B}$ given by 
\[\Delta(\lambda_0\otimes \cdots\otimes \lambda_{n+1})=\sum_{i=0}^n
(\lambda_0\otimes \cdots\otimes \lambda_i\otimes 1)\otimes (1\otimes
\lambda_{i+1}\otimes \cdots\otimes \lambda_{n+1})\] (for example see
\cite[1.2]{S}).  If $\eta$ and $\theta$ in $\HH^n(\L)$ and
$\HH^m(\L)$, respectively, are represented by $\eta\colon B^n\to \L$
and $\theta\colon B^m\to \L$, then the cup product $\eta\cup \theta$
in $\HH^{n+m}(\L)$ is given by the following composition of maps
\[\mathbb{B}\extto{\Delta} \mathbb{B}\otimes_\L
\mathbb{B}\extto{\eta\otimes \theta}\L\otimes_\L \L\extto{\nu}\L,\] 
where $\nu\colon \L\otimes_\L \L\to \L$ is the multiplication map.
We see that the cup product is
\[\eta\cup \theta((\lambda_0\otimes \cdots\otimes \lambda_{n+m+1})) =
\eta(\lambda_0\otimes \cdots\otimes \lambda_n\otimes 1)\theta(1\otimes
\lambda_{n+1}\otimes \cdots\otimes \lambda_{n+m+1}).\]

In \cite{SW} it is shown that any projective $\L^e$-resolution
$\mathbb{X}$ of $\L$ gives rise to a ``cup product'', which coincides
with the ordinary cup product. Let $\mathbb{X}$ be a projective
$\L^e$-resolution of $\L$. There exists a chain map $D\colon
\mathbb{X}\to \mathbb{X}\otimes_\L\mathbb{X}$ lifting the identity,
which is unique up to homotopy. Siegel and Witherspoon define a cup
product of two elements $\eta$ in $\HH^n(\L)$ and $\theta$ in
$\HH^m(\L)$ as above using the composition
\[\mathbb{X}\extto{D}\mathbb{X}\otimes_\L\mathbb{X}
\extto{\eta\otimes\theta} 
\L\otimes_\L \L\extto{\nu}\L,\]
and note that it is independent of the projective resolution
$\mathbb{X}$ of $\L$ and the chain map $D$. In Section \ref{section:2}
we give an explicit formula for $D$ for the minimal projective
$\L^e$-resolution of $\L$ constructed here. 

In \cite{Ge} it is shown that the cup product and the Yoneda product
coincide. We give a proof of this below. Let $\eta\ast \theta$ denote
the Yoneda product of $\eta$ and $\theta$ viewing $\HH^*(\L)$ as
$\Ext^*_{\L^e}(\L,\L)$.

\begin{prop}\label{prop:product}
Let $\eta$ and $\theta$ be in $\HH^n(\L)$ and $\HH^m(\L)$,
respectively, and suppose that they are represented by $\eta\colon
B^n\to \L$ and $\theta\colon B^m\to \L$. Then \[\eta\ast \theta=
\eta\cup \theta.\] In particular, the cup product and the Yoneda
product on $\HH^*(\L)$ coincide.
\end{prop}
\begin{proof}
Let $\eta$ and $\theta$ be in $\HH^n(\L)$ and $\HH^m(\L)$
respectively. Then $\eta$ and $\theta$ can be represented as
$\L^e$-maps $\eta\colon B^n\to \L$ and $\theta\colon B^m\to \L$. The
Yoneda product $\eta\ast\theta$ is given via constructing a lifting of
$\theta$ to a chain map $\widetilde{\theta}=\{\theta^i\}_i\colon
\mathbb{B}\to \mathbb{B}[n]$, where $\mathbb{B}[n]$ denotes the degree
$n$ shift of the complex $\mathbb{B}$. Define $\widetilde{\theta}$ to
be the composition
\[\mathbb{B}\extto{\Delta} \mathbb{B}\otimes_\L
\mathbb{B}\extto{\mathbb{B}\otimes\theta}
\mathbb{B}\otimes_\L\L[n]\extto{\nu}\mathbb{B}[n]\] 
It is clear from the definition of $\widetilde{\theta}$ that it is a
morphism of complexes. Furthermore, direct computations shows that
$\widetilde{\theta}$ is a lifting of $\theta$ to a chain map. Immediately
from the definition of $\widetilde{\theta}$ we infer that 
\begin{alignat}{2}
\eta\ast\theta & = \eta\widetilde{\theta}
               & & = \eta\nu(\mathbb{B}\otimes\theta)\Delta\notag\\
& = \nu(\eta\otimes\L[n])(\mathbb{B}\otimes\theta)\Delta 
& & = \nu(\eta\otimes\theta)\Delta\notag\\ 
& = \eta\cup \theta. & & \notag
\end{alignat}
\end{proof}

\section{The multiplicative structure of the Hochschild cohomology
ring}\label{section:2}

In this section the multiplicative structure of the Hochschild
cohomology ring of a Koszul algebra, $\L=kQ/I$, is found using the
comultiplicative structure of a minimal projective resolution of
$\L_0$ as a right $\L$-module as given in equation \eqref{eq:comult}
in the introduction. More precisely, given two homogeneous elements in
the Hochschild cohomology ring, we give a closed formula for their
product as a map from the appropriate projective module in a
projective resolution of $\L$ over $\L^e$ to $\L$. Note that this map
may be non-zero, but represent zero in Hochschild cohomology, since
the residue class of the map in Hochschild cohomology still needs to
be computed (see Example \ref{ex:1}). The crucial ingredient in the
proof is the minimal projective resolution of $\L$ as a right
$\L^e$-module constructed in \cite{GHMS} and the inclusion of this
resolution into the reduced bar resolution.

We begin by recalling definitions, notation, and results that we need
in this section. In this section $\L=kQ/I$ denotes a Koszul algebra
over some field $k$, where $Q$ is a finite quiver and $I\subseteq
J^2$. Here $J$ denotes the ideal in $kQ$ generated by the arrows, and
let $\L_0=\L/(J/I)$. If $M$ is a $kQ$-module and $m$ is in $M$, then
$\overline{m}$ denotes the natural residue class of $m$ in $M/MI$.

In \cite{GHMS}, it was shown that a minimal projective resolution
$(\mathbb{P},\delta)$ of $\L$ as a right $\L^e$-module can be
constructed from knowledge of a minimal projective resolution of
$\L_0$ over $\L$ as a right $\L$-module.  Since we need the results
from \cite{GHMS}, we summarize them below.

A minimal projective resolution $(\mathbb{F},d)$ of $\L_0$ as a right
$\L$-module can be constructed from a sequence of right ideals in
$R=kQ$ as follows: There exists a choice of integers $t_n$ in $\Z$ and
of uniform elements $\{f^n_i\}_{i=0}^{t_n}$ in $R$, for all $n\geq 0$,
with
\[\cdots\subseteq \amalg_{i=0}^{t_n} f^n_iR \subseteq
\amalg_{i=0}^{t_{n-1}}  
f^{n-1}_iR\subseteq \cdots \subseteq \amalg_{i=0}^{t_1}
f^1_iR\subseteq \amalg_{i=0}^{t_0}f^0_iR=R\] such that
$F^n=\amalg_{i=0}^n f^n_iR/\amalg_{i=0}^n f^n_iI$ for all $n\geq 0$
and $d^n\colon F^n\to F^{n-1}$ is induced by the inclusion
$\amalg_{i=0}^n f^n_iR\hookrightarrow \amalg_{i=0}^{n-1} f^{n-1}_iR$
(see \cite{GSZ}). We choose $\{f^0_i\}_{i=0}^{t_0}$ to be the set of
the vertices in $Q$ and $\{f^1_i\}_{i=0}^{t_1}$ to be the set of the
arrows in $Q$. Moreover $\{f^2_i\}_{i=0}^{t_2}$ is a minimal set of
homogeneous generators of degree two for $I$.

When we have a minimal projective resolution $(\mathbb{F},d)$ of
$\L_0$ as a right $\L$-module obtained from a set $\{f^n_i\}$, we
simply say that the set $\{f^n_i\}$ \emph{defines a minimal projective 
resolution}.

It is shown in \cite{GHMS} that a minimal projective resolution
$(\mathbb{P},\delta)$ of $\L$ over $\L^e$ is given by
$P^n=\amalg_{i=0}^{t_n} \L\frako(f^n_i)\otimes_k \frakt(f^n_i)\L$ for
$n\geq 0$, where the differential $\delta^n\colon P^n\to P^{n-1}$ is
now described for $n\geq 1$: In \cite{GHMS} it was shown that the sets
$\{f^n_i\}$ have a ``comultiplicative structure''; namely, there are
elements $c_{pq}(n,i,r)$ in $k$ such that
\begin{equation}
f^n_i  = \sum_{p=0}^{t_r}\sum_{q=0}^{t_{n-r}}
c_{pq}(n,i,r)f^r_{p}f^{n-r}_{q}\tag{1}
\end{equation}
for all $n\geq 1$, all $i$ in $\{0,1,\ldots,t_n\}$ and all $r$ in
$\{0,1,\ldots,n\}$.  To simplify notation, we sometimes write
$\sum_{p,q}$ for $\sum_{p=0}^s\sum_{q=0}^{s'}$ when the bounds of $p$
and $q$ are clear. For $n\geq 0$ and $i$ with $0\leq i\leq t_n$
define $\varepsilon^n_i=(0,\ldots,0,\underbrace{\frako(f^n_i)\otimes_k
\frakt(f^n_i)}_{\text{$i$-th component}},0,\ldots,0)$ in $P^n$.  Then
the differential $\delta^n\colon P^n\to P^{n-1}$ is given by
\[\delta^n(\varepsilon^n_i)=\sum_{j=0}^{t_{n-1}}\left(\sum_{p=0}^{t_1} 
c_{pj}(n,i,1)\overline{f^1_{p}}\varepsilon^{n-1}_j +(-1)^n
\sum_{q=0}^{t_1} 
c_{jq}(n,i,n-1)\varepsilon^{n-1}_j\overline{f^1_{q}}\right)\]
for $i=0,1,\ldots,t_n$ and $n\geq 1$, and $\delta^0\colon
\amalg_{i=0}^{t_0}\L e_i\otimes_k e_i\Lambda \to \L$ is the
multiplication map. Note that $\overline{f^n_j}$ denotes the residue
class of $f^n_j$ in
$\amalg_{i=0}^{t_n}f^n_iR/\amalg_{i=0}^{t_n}f^n_iI$. 

Recall that the reduced bar resolution $(\mathbb{B},d)$ of $\L$ over
$\L^e$ from the previous section is given by
$B^n=\L^{\otimes_{\L_0}(n+2)}$ and the differential $d^n\colon B^n\to
B^{n-1}$ is given by
\begin{multline}
d^n(\lambda_0\otimes\lambda_1\otimes\cdots\otimes \lambda_{n+1})
=\notag\\ \sum_{i=0}^n
(-1)^i(\lambda_0\otimes\lambda_1\otimes\cdots\otimes\lambda_{i-1}
\otimes \lambda_i\lambda_{i+1}\otimes \lambda_{i+2}\otimes\cdots
\otimes \lambda_{n+1}).\notag
\end{multline}
Define $\partial_i\colon B^n\to B^{n-1}$ by 
\[\partial_i(\lambda_0\otimes \cdots\otimes
\lambda_{n+1})=\lambda_0\otimes  
\cdots\otimes \lambda_{i-1}\otimes \lambda_i\lambda_{i+1}\otimes
\lambda_{i+2}\otimes \cdots\otimes \lambda_{n+1}\] for
$i=0,1,\ldots,n$. Then $d^n=\sum_{i=0}^n (-1)^i\partial_i$.

Identify $\L_0$ with the subalgebra of $\L$ generated by the vertices.
View $R$ as the tensor algebra $T_{\L_0}(V)$, where $V$ is the
$\L_0$-bimodule generated by the arrows in $Q$. Since $\L$ is a Koszul
algebra, each $f^n_i$ is a linear combination of paths of length
$n$. Hence each $f^n_i$ can be viewed uniquely as an element in
$V^{\otimes_{\L_0}n}$ for all $n$ and $i$; see the definition of
$\psi(f^n_i)$ in the proposition below. Now we show that this enables
us to find a natural embedding of $\mathbb{P}$ as a subcomplex of
$\mathbb{B}$.

\begin{prop}
Define $\mu_n\colon P^n\to
B^n$ by
\[\mu_n(\varepsilon^n_i) = 1\otimes_{\L_0} \psi(f^n_i)\otimes_{\L_0} 1\] 
where
\[\psi(f^n_i)=\sum c_{ij_1j_2\ldots j_n}\overline{f^1_{j_1}}\otimes
\cdots \otimes \overline{f^1_{j_n}}\]
when $f^n_i=\sum  c_{ij_1j_2\ldots j_n}f^1_{j_1}\cdots f^1_{j_n}$. 

Then $\{\mu_n\}_{n\geq 0}$ is a chain map from $(\mathbb{P},\delta)$
to $(\mathbb{B},d)$.
\end{prop}
\begin{proof}
If $f^n_i=\sum_{p,q} c_{pq}(n,i,r)f^r_{p}f^{n-r}_{q}$, then 
$\psi(f^n_i)=\sum_{p,q} c_{pq}(n,i,r) \psi(f^r_{p})\otimes
\psi(f^{n-r}_{q})$. 

We have that 
\begin{align}
f^n_i & = \sum_{p,q} c_{pq}(n,i,r)f^r_{p}f^{n-r}_{q}\notag\\ & =
      \sum_{p,q} c_{pq}(n,i,r)f^r_{p}
      (\sum_{p',q'}c_{p'q'}(n-r,q,2)f^2_{p'}f^{n-r-2}_{q'})\notag
\end{align}
for any $r$ with $0\leq r\leq n$. It follows from this that 
\[\partial_r\mu_n(\varepsilon^n_i) = 
\sum_{p,q} c_{pq}(n,i,r)\otimes \psi(f^r_{p})\otimes 
      (\sum_{p',q'}c_{p'q'}(n-r,q,2) \overline{f^2_{p'}}\otimes
\psi(f^{n-r-2}_{q'}))\otimes 1\]
whenever $n\geq 2$ and $0\leq r\leq n-2$. We see that 
$\partial_r\mu_n(\varepsilon^n_i)$  is zero since
$f^2_{i_{p'}}$ is in the ideal $I$ for all $p'$. We infer from this that
\begin{multline}
d^n\mu_n(\varepsilon^n_i) = \sum_{p,q}
c_{pq}(n,i,1)\overline{f^1_{p}}\otimes\psi(f^{n-1}_{q})\otimes
1\\ \quad + (-1)^n \sum_{p',q'} c_{p'q'}(n,i,n-1)\otimes 
\psi(f^{n-1}_{p'})\otimes \overline{f^1_{q'}}\notag
\end{multline}
and consequently that $d^n\mu_n=\mu_{n-1}\delta^n$. This shows that
the chain map $\{\mu_n\}_{n\geq 0}$ defines $(\mathbb{P},\delta)$ as a
subcomplex of the reduced bar resolution $(\mathbb{B},d)$ from the
previous section.
\end{proof}

A projective resolution of $\L$ over $\L^e$ different from the bar
resolution was first described for a Koszul algebra $\L$ in
\cite[Section 3.7]{P} (choose $R=L=A$ in $K_*(R,A,L)$ in the notation
of \cite{P}). Furthermore, an embedding of this resolution into the
bar resolution is given in \cite[Proposition 3.9]{P}.

Given two homogeneous elements $\eta$ and $\theta$ in the cohomology
of $\Hom_{\L^e}(\mathbb{B},\L)$ their product $\eta\ast\theta$ in
$\HH^*(\L)$ is given as the composition of the maps
\[\mathbb{B}\extto{\Delta} \mathbb{B}\otimes_\L\mathbb{B}
\extto{\eta\otimes \theta} \L\otimes_\L\L\extto{\nu}\L.\] 
Let $\mu\colon \mathbb{P}\to
\mathbb{B}$ be the inclusion found above, and let $\pi\colon
\mathbb{B}\to \mathbb{P}$ be a chain map such that
$\pi\mu=\id_{\mathbb{P}}$. 

Suppose that the image $\Delta(\mu\mathbb{P})$ is contained in
$\mu\mathbb{P}\otimes_\L\mu\mathbb{P}\subseteq
\mathbb{B}\otimes_\L\mathbb{B}$.  In that case, $\Delta$ induces a map
$\mathbb{P}\to \mathbb{P}\otimes_\L \mathbb{P}$ that we denote by
$\Delta'$. In particular $\Delta\mu=(\mu\otimes\mu)\Delta'$. By the
results from \cite{SW} mentioned prior to Proposition
\ref{prop:product}, the cup product defined by $\Delta'$ gives the
multiplication in $\HH^*(\L)$.

Next we show that this is indeed the case; that is,
$\Delta(\mu\mathbb{P})$ is contained in
$\mu\mathbb{P}\otimes_\L\mu\mathbb{P}$. To this end define
$\Delta'\colon \mathbb{P}\to \mathbb{P}\otimes_\L \mathbb{P}$ by
\[\Delta'(\varepsilon^n_i) =
\sum_{r=0}^n \sum_{p=0}^{t_r}\sum_{q=0}^{t_{n-r}} c_{pq}(n,i,r) 
\varepsilon^r_p\otimes_\L \varepsilon^{n-r}_q.\] 
\begin{prop}
Let $\Delta$ and $\Delta'$ be as above. Then
\[\Delta\mu=(\mu\otimes\mu)\Delta'\]
\end{prop}
\begin{proof}
Let $s_r\colon \mathbb{B}\to \mathbb{B}\otimes_\L\mathbb{B}$ be
given by 
\[s_r(\lambda_0\otimes\cdots\otimes\lambda_{n+1}) =
(\lambda_0\otimes\cdots\otimes\lambda_r\otimes 1)\otimes_\L (1\otimes
\lambda_{r+1}\otimes\cdots\otimes\lambda_{n+1})\]
for $r=0,1,\ldots,n$. Then we have that $\Delta=\sum_{r=0}^n
s_r$.  Furthermore, 
\begin{align}
s_r(\mu_n(\varepsilon^n_i)) & =
\sum_{p,q}(c_{pq}(n,i,r)\otimes\psi(f^r_{p})\otimes
1) \otimes_\L (1\otimes \psi(f^{n-r}_{q})\otimes 1)\notag\\
& = \sum_{p,q}(c_{pq}(n,i,r)\mu_r(\varepsilon^r_p)
\otimes_\L \mu_{n-r}(\varepsilon^{n-r}_q)\notag\\
& = (\mu_r\otimes\mu_{n-r})\Delta'(\varepsilon^n_i).\notag
\end{align}
Since $\Delta =\sum_{r=0}^n s_r$, we have that
$\Delta\mu_n(\varepsilon^n_i) = (\mu\otimes\mu)\Delta'(\varepsilon^n_i)$.
Hence the claim follows. 
\end{proof}
Combining the results above we obtain our result on the multiplicative
structure for the Hochschild cohomology ring of a Koszul algebra.  

\begin{thm}\label{thm:multformula}
Let $\L=kQ/I$ with $I\subseteq J^2$ be a Koszul algebra over a field
$k$, where $Q$ is a finite quiver.  Suppose that $\eta\colon P^n\to
\L$ and $\theta\colon P^m\to \L$ represent elements in $\HH^*(\L)$ and
that they are given by $\eta(\varepsilon^n_i)=\lambda_i$ for
$i=0,1,\ldots,t_n$ and $\theta(\varepsilon^m_i)=\lambda_i'$ for
$i=0,1,\ldots,t_m$. Then
\[(\eta\ast \theta)(\varepsilon^{m+n}_i)
= \sum_{p=0}^{t_n}\sum_{q=0}^{t_m}c_{pq}(n+m,i,n)
\lambda_{p}\lambda_{q}'\]
for all $i=0,1,\ldots, t_{m+n}$. 
\end{thm}
\begin{proof}
By Proposition \ref{prop:product}, we have that $\eta\ast \theta =
\eta\cup \theta$. Using the results from \cite{SW} reviewed prior to
Proposition \ref{prop:product} we infer that 
\begin{align}
\eta\ast \theta(\varepsilon^{m+n}_i) & = 
\nu(\eta\otimes \theta)\Delta'(\varepsilon^{m+n}_i)\notag\\
& = \nu(\eta\otimes \theta)(\sum_{r=0}^{n+m}
\sum_{p,q}c_{pq}(n+m,i,r)\varepsilon^r_p 
\otimes_\L \varepsilon^{n-r}_q)\notag\\ 
& = \sum_{p,q}c_{pq}(n+m,i,n) \eta(\varepsilon^n_p)
\theta(\varepsilon^m_q)\notag\\
& = \sum_{p,q}c_{pq}(n+m,i,n) \lambda_p\lambda_q'.\notag
\end{align}
This completes the proof of the theorem.
\end{proof}

\section{Structure constants for the Koszul dual}\label{section:3}

In this section we give a second application of the comultiplication
formula \eqref{eq:comult}. We show that, for a Koszul algebra
$\L=kQ/I$, equation \eqref{eq:comult} gives the structure constants
for a basis for the Koszul dual of $\L$. Of course the Koszul dual can
immediately be given by generators and relations. But for
computational purposes, it is sometimes important to find a $k$-basis
and the structure constants for that basis.

Suppose that the sets $\{f^n_i\}$ define a minimal projective
resolution $(\mathbb{F},e)$ of $\L/\rrad$ as a right $\L$-module. Note
that $F^n=\amalg_{i=0}^{t_n}f^nR/\amalg_{i=0}^{t_n}f^nI$. Define
$\widehat{f}^n_i\colon F^n\to \L/\rrad$ for $i=0,1,\ldots, t_n$ by
$\widehat{f}^n_i(\overline{f^n_j})=\delta_{ij}\frakt(f^n_i)$ for all
$j=0,1,\ldots,t_n$. Since $(\mathbb{F},e)$ is a minimal projective
resolution, it is immediate that the elements
$\{\widehat{f}^n_i\}_{i=0}^{t_n}$ form a basis for
$\Ext^n_\L(\L/\rrad,\L/\rrad)$ as a vector space over $k$. The goal of
this section is to find the structure constants for this particular
$k$-basis for the Koszul dual of $\L$. We note that $\L$ is given by
generators and relations, whereas the multiplication in the Koszul
dual is given by structure constants. Since $\L$ is the Koszul dual of
its Koszul dual, we may construct the set $\{ f^n_i\}$ for the Koszul
dual and then the results here can be applied to find a $k$-basis of
$\L$ and its structure constants.

\begin{thm}\label{thm:structconst}
Let $\L=kQ/I$ be a Koszul algebra over some field $k$, where $Q$ is a
finite quiver and $I\subseteq J^2$. Denote by $\{f^n_i\}_{i=0}^{t_n}$
elements in $kQ$ defining a minimal projective resolution
$(\mathbb{F},e)$ of $\L/\rrad$ as a right module over $\L$.  Let the
elements $\widehat{f}^n_i\colon F^n\to \L/\rrad$ for $i=0,1,\ldots,
t_n$ be given by
$\widehat{f}^n_i(\overline{f^n_j})=\delta_{ij}\frakt(f^n_i)$, for all
$j=0,1,\ldots,t_n$. The set $\{\widehat{f}^n_i\}_{i=0}^{t_n}$
represents a $k$-basis of $\Ext^n_\L(\L/\rrad,\L/\rrad)$.

Then the product $\widehat{f}^m_i\widehat{f}^n_j$ in
$\Ext^*_\L(\L/\rrad,\L/\rrad)$ of the two elements $\widehat{f}^m_i$
and $\widehat{f}^n_j$ in degree $m$ and $n$ in
$\Ext^*_\L(\L/\rrad,\L/\rrad)$, respectively, is given by
\[\widehat{f}^m_i\widehat{f}^n_j=\sum_{l=0}^{t_{m+n}}
c_{ji}(m+n,l,n)\widehat{f}^{m+n}_l.\] 
\end{thm}
\begin{proof}
A lifting of the map $\widehat{f}^n_j\colon F^n\to \L/\rrad$ to a map
$g_0\colon F^n\to F^0$ is given by
$g_0(\overline{f^n_t})=\delta_{jt}\frakt(f^n_j)$. We
want to find a lifting of $g_0$ to a chain map
$\mathbb{F}[n]\to \mathbb{F}$, where $\mathbb{F}[i]$ denotes the
$i$-th shift of the complex $\mathbb{F}$. Define
$g_r\colon F^{n+r}\to F^r$ by letting
\[g_r(\overline{f^{n+r}_i})=\sum_{q=0}^{t_r}\overline{f^r_q}
c_{jq}(n+r,i,n)\]
for $i=0,1,\ldots,t_{n+r}$. Then it is clear that
$g_0e^{n+1} = e^1 g_1$.  To show that
$g_re^{n+r+1} = e^{r+1}g_{r+1}$ for $r\geq
1$ consider the following equalities:
\begin{align}
f^{n+r+1}_l & = \sum_{p=0}^{t_{n+r}}\sum_{q=0}^{t_1} 
f^{n+r}_pf^1_q c_{pq}(n+r+1,l,n+r)\notag\\
 & = \sum_{p=0}^{t_{n+r}}\sum_{q=0}^{t_1}\sum_{x=0}^{t_n}\sum_{y=0}^{t_r}
f^n_x f^r_y f^1_qc_{xy}(n+r,p,n)c_{pq}(n+r+1,l,n+r)\notag 
\end{align}
and
\begin{align}
f^{n+r+1}_l & = \sum_{x=0}^{t_n}\sum_{v=0}^{t_{r+1}} f^n_x f^{r+1}_v
c_{xv}(n+r+1,l,n) \notag\\
& = \sum_{x=0}^{t_n}\sum_{v=0}^{t_{r+1}}\sum_{y=0}^{t_r}\sum_{q=0}^{t_1}
f^n_x f^r_y f^1_q c_{xv}(n+r+1,l,n)c_{yq}(r+1,v,r).\notag
\end{align}
Since $\sum_{x=0}^{t_n} f^n_xR$ is a direct sum, 
for $x=j$, we have that 
\begin{multline}
\sum_{p=0}^{t_{n+r}}\sum_{q=0}^{t_1}\sum_{y=0}^{t_r} f^r_y
f^1_qc_{jy}(n+r,p,n)c_{pq}(n+r+1,l,n+r) = \\
\sum_{v=0}^{t_{r+1}}\sum_{y=0}^{t_r}\sum_{q=0}^{t_1} f^r_y f^1_q
c_{jv}(n+r+1,l,n)c_{yq}(r+1,v,r) \notag
\end{multline}
as elements in $\amalg_{y=0}^{t_r}f^r_yR$. Taking these equalities
in $\amalg_{y=0}^{t_r}f^r_yR/\amalg_{y=0}^{t_r}f^r_yI$, the left hand
side of this equality is $g_re^{n+r+1}(\overline{f^{n+r+1}_l})$ and
the right hand side is $e^{r+1}g_{r+1}(\overline{f^{n+r+1}_l})$ for
all $l=0,1,\ldots,t_{n+r+1}$, hence $g_re^{n+r+1} =
e^{r+1}g_{r+1}$. We have that
$\widehat{f}^m_i\widehat{f}^n_j=\widehat{f}^m_ig_m$, so that
$\widehat{f}^m_i\widehat{f}^n_j= \sum_{l=0}^{t_{m+n}}
c_{ji}(m+n,l,n)\widehat{f}^{m+n}_l$. This completes the proof.
\end{proof}

\section{The graded centre and Hochschild cohomology}\label{section:4}

For a finite dimensional algebra or graded algebra $\L=kQ/I$ the image
of the natural map $\varphi_{\L/\rrad}\colon \HH^*(\L)\to
\Ext^*_\L(\L/\rrad,\L/\rrad)=E(\L)$ is shown to be contained in the
graded centre $Z_\gr(E(\L))$ in \cite{SS}. Here $\varphi_{\L/\rrad}$
is induced from $\L/\rrad\otimes_\L-\colon (\mathbb{P},\delta)\to
(\L/\rrad\otimes_\L \mathbb{P},1\otimes \delta)$.  It was
independently observed by Buchweitz and Green-Snashall-Solberg that
when $\L$ is a Koszul algebra, the image is in fact equal to
$Z_\gr(E(\L))$. This was obtained by Buchweitz as a part of a more
general isomorphism between the Hochschild cohomology ring of $\L$ and
the graded Hochschild cohomology ring of the Koszul dual. This
isomorphism has since been generalized by Keller \cite{K}. We note
that the image of $\varphi_{\L/\rrad}$ is in general strictly
contained in $Z_\gr(E(\L))$, see \cite[Example 7.6]{GSS}. In this
section we give an elementary proof of the fact that the image is
$Z_\gr(E(\L))$ for a Koszul algebra $\L$ using the results from the
previous sections.

\begin{thm}
Let $\L=kQ/I$ be a Koszul algebra. The image of the natural map
$\varphi_{\L/\rrad}\colon \HH^*(\L)\to E(\L)$ is the graded centre
$Z_\gr(E(\L))$. 
\end{thm}
\begin{proof}
We keep the notation developed in the previous sections.  Let
$z=\sum_{i=0}^{t_n}\alpha_i \widehat{f}^n_i$ with $\alpha_i$ in
$k$. Since $\L$ is a Koszul algebra, $E(\L)$ is generated in degrees
$0$ and $1$. Hence $z$ is in $Z_\gr(E(\L))$ if and only if (A)
$\widehat{f}^1_jz=(-1)^nz\widehat{f}^1_j$ for all $j=0,1,\ldots,t_1$
and (B) $\widehat{f}^0_lz=z\widehat{f}^0_l$ for all
$l=0,1,\ldots,t_0$. Using the structure constants found in Theorem
\ref{thm:structconst}, we have that
\[\widehat{f}^1_jz = \sum_{i=0}^{t_n} \alpha_i
\widehat{f}^1_j\widehat{f}^n_i = \sum_{i=0}^{t_n} \alpha_i
\sum_{l=0}^{t_{n+1}} c_{ij}(n+1,l,n) \widehat{f}^{n+1}_l\] and
\[z\widehat{f}^1_j =
\sum_{i=0}^{t_n}\alpha_i\widehat{f}^n_i\widehat{f}^1_j =
\sum_{i=0}^{t_n}\alpha_i\sum_{l=0}^{t_{n+1}}
c_{ji}(n+1,l,1)\widehat{f}^{n+1}_l.\]
Suppose that $z$ is in $Z_\gr(E(\L))$. Then the condition (B) above
implies that $\frako(f^n_i)=\frakt(f^n_i)$ for all $i$ such that
$\alpha_i\neq 0$ and 
\[\sum_{i=0}^{t_n}\alpha_ic_{ij}(n+1,l,n)=(-1)^n\sum_{i=0}^{t_n}\alpha_i
c_{ji}(n+1,l,1).\] 
Define $\eta\colon P^n\to \L$ by letting
$\eta(\varepsilon^n_i)=\alpha_i\frakt(f^n_i)$ for
$i=0,1,\ldots,t_n$. Then 
\begin{align}
\eta\delta^{n+1}(\varepsilon^{n+1}_l) & =
\sum_{i=0}^{t_n}\sum_{j=0}^{t_1}
c_{ji}(n+1,l,1)\overline{f^1_j}\alpha_i+
(-1)^{n+1}\sum_{i=0}^{t_n}\sum_{q=0}^{t_1}
c_{iq}(n+1,l,n)\overline{f^1_q}\alpha_i\notag\\ 
& = \sum_{j=0}^{t_1}\sum_{i=0}^{t_n}(c_{ji}(n+1,l,1) + (-1)^{n+1}
c_{ij}(n+1,l,n))\alpha_i\overline{f^1_j}\notag\\ & = 0.\notag
\end{align}
for all $l=0,1,\ldots, t_{n+1}$ and $j=0,1,\ldots,t_1$. 
The last equality follows since $z$ is in $Z_\gr(E(\L))$. Hence $\eta$
is in $\HH^n(\L)$, and $\varphi_{\L/\rrad}(\eta)=z$. This completes
the proof.
\end{proof}
In \cite{SS} it is conjectured that, if $\L$ is a finite dimensional
algebra over a field $k$, then the Hochschild cohomology ring of $\L$
modulo the ideal generated by the homogeneous nilpotent elements is a
finitely generated commutative algebra over $k$. As a consequence of
the previous result, for a finite dimensional Koszul algebra $\L$, the
conjecture is equivalent to the conjecture that $Z_\gr(E(\L))$, modulo
the ideal generated by homogeneous nilpotent elements, is a finitely
generated algebra over $k$. This can be seen by noting that the kernel
of $\varphi_{\L/\rrad}$ is contained in the ideal generated by
homogeneous nilpotent elements \cite{SS}. We apply these ideas in
Example \ref{ex:5.3}.

\section{Examples}

\begin{example}\label{ex:1}
This example shows that the multiplication formula in Theorem
\ref{thm:multformula} only gives a representative of a product in
Hochschild cohomology as a map from the appropriate projective to the
Koszul algebra.  Thus, for example, to find if two products are equal,
the residue classes must be computed.
 
Let $\L=k\langle x,y\rangle/(x^2,xy+yx)$ for a field $k$ of
characteristic different from $2$. Since $\{x^2,xy+yx\}$ is a
quadratic Gr\"obner basis for the ideal generated by the relations
under the length lexicographic order with $x>y>1$, the algebra $\L$ is
Koszul \cite{GH}. In this example $t_0=0$, and $t_n=1$ for all $n\geq
1$. We have that $f^0_0=1$, and $f^n_0=x^n$ and
$f^n_1=\sum_{\substack{a,b\geq 0\\ a+b=n-1}} x^ayx^b$ for $n\geq
1$. Set $f^0_1=0$. For $n\geq 1$ and $r$ with $0\leq r\leq n$, we see
that $f^n_0=f^r_0f^{n-r}_0$ and that
$f^n_1=f^r_0f^{n-r}_1+f^r_1f^{n-r}_0$. Hence,
$c_{00}(n,0,r)=c_{01}(n,1,r)=c_{10}(n,1,r)=1$ and all other
$c_{pq}(n,i,r)=0$. In particular, $c_{01}(2,1,1)=c_{10}(2,1,1)=1$.
Let $(\mathbb{P},\delta)$ be the minimal projective resolution of $\L$
as a right $\L^e$-module described in Section \ref{section:2}.  Then
$P^1=\varepsilon^1_0\L^e\amalg \varepsilon^1_1\L^e$.  Define $\eta\colon
P^1\to \L$ by $\eta(\varepsilon^1_i)=\begin{cases} xy, & i=0\\ y, & i=1
\end{cases}$, and $\theta\colon P^1\to \L$ by
$\eta(\varepsilon^1_i)=\begin{cases}
0, & i=0\\
y, & i=1
\end{cases}$. The reader may check that both $\eta$ and $\theta$
represent non-zero elements in Hochschild cohomology. On the other
hand, the map $\eta\ast\theta$ is non-zero, but represents zero in
$\HH^2(\L)$.
\end{example}

\begin{example}
Let $Q$ be a finite quiver, which is not an oriented cycle. Recall
that $J$ denotes the ideal of $kQ$ generated by the arrows of $Q$. Let
$\L=kQ/J^2$, which is a Koszul algebra. Then $t_n+1$ is the number of
paths of length $n$ for all $n\geq 0$, and we choose the set
$\{f^n_i\}_{i=0}^{t_n}$ to be the set of all paths of length $n$ in
some order. Given $n$ and $r$ with $0\leq r\leq n$, and $i$ with
$0\leq i\leq t_n$, there exist unique $j$ and $l$ with $0\leq j\leq
t_r$ and $0\leq l \leq t_{n-r}$ such that
$f^n_i=f^r_jf^{n-r}_l$. Hence $c_{pq}(n,i,r)=\begin{cases} 1, &
\text{if $p=j$ and $q=l$}\\ 0, & \text{otherwise}\end{cases}.$ A
classic result states that the Koszul dual of $\L$ is $kQ^\op$, where
$Q^\op$ is the opposite quiver of $Q$. If we write
$f^n_i=f^1_{j_1}\cdots f^1_{j_n}$ as a product of arrows, then by
Theorem \ref{thm:structconst} the element
$\widehat{f}^n_i=\widehat{f}^1_{j_n}\cdots \widehat{f}^1_{j_1}$ and
therefore the Koszul dual of $\L$ is isomorphic to $kQ^\op$.

For a discussion of the Hochschild cohomology ring of $\L$ we refer to
\cite[Lemma 7.4]{GSS} and the remarks following it or to \cite{C}.  
\end{example}

\begin{example}\label{ex:5.3}
Let $\L=k\langle x,y,z\rangle/(x^2,y^2,z^2,xy+ayx,xz+bzx,yz+czy)$ be
the quantum exterior algebra over a field $k$ with $a$, $b$ and $c$
non-zero elements in $k$. This algebra is easily seen to be
selfinjective of dimension $8$ over $k$. In the length lexicographic
order with $x>y>z>1$, the set $\{x^2,y^2,z^2,xy+ayx,xz+bzx,yz+czy\}$
is a quadratic Gr\"obner basis for the ideal it generates, and hence
$\L$ is a Koszul algebra \cite{GH}. The Koszul dual of $\L$ is quantum
$3$-space, $R=k\langle x,y,x\rangle/(yx-axy,zx-bxz,zy-cyz)$. Note that
$R$ is Artin-Schelter regular of global dimension $3$.

The Hochschild cohomology ring of the algebra $\Sigma=k\langle
x,y\rangle/(x^2,xy+qyx,y^2)$ for $q$ in $k$ was studied in detail
in \cite{BGMS}. In that paper it was shown that when $q$ is not a root
of unity, then the Hochschild cohomology groups $\HH^n(\Sigma)$
vanish for $n\geq 3$. In this example we show that a similar
phenomena occurs for $\L$. In particular, for certain values of $a$,
$b$, and $c$, we show that $\HH^n(\L)$ vanishes for $n\geq 4$.

For ease of notation, we triply index the set $\{f^n_i\}$ as the set
$\{f^n_{u,v,w}\}$, where $u+v+w=n$. A choice of the elements $\{
f^n_{u,v,w}\}$ defining a minimal resolution of $\L/\rrad$ as a right
$\L$-module can be given inductively as follows: $f^0_{0,0,0}=1$ and for
non-negative integers $u$, $v$, and $w$ with $u+v+w=n$,
\[f^n_{u,v,w}=a^vb^wf^{n-1}_{u-1,v,w}x+c^wf^{n-1}_{u,v-1,w}y+ 
f^{n-1}_{u,v,w-1}z\] with the requirement that
$f^n_{-1,v,w}=f^n_{u,-1,w}=f^n_{u,v,-1}=0$ and that $f^n_{u,v,w}=0$
for $u+v+w>n$. The reader may check that for each $n$, the number of
$f^n_{u,v,w}$'s is $\binom{n+2}{n}$. 

Let $(\mathbb{P},\delta)$ be the minimal projective $\L^e$-resolution
of $\L$ defined in Section \ref{section:2}.  To compute the Hochschild
cohomology groups, we compute the dimension of the image of the map
\[(\delta^n)^*=\Hom_{\L^e}(\delta^n,\L)\colon \Hom_{\L^e}(P^{n-1},\L)\to
\Hom_{\L^e}(P^n,\L).\] One finds that the dimension of
$\Im(\delta^n)^*$ is given by $2n^2+4n+1$ when $n\geq 3$ if, for
example, $a$, $b$ and $c$ are algebraically independent elements of
$k$. This implies that the Hochschild cohomology groups $\HH^n(\L)$
vanish for $n\geq 4$, when $a$, $b$ and $c$ are algebraically
independent elements of $k$. 

In the case when $a$, $b$ and $c$ are algebraically independent, it
is easy to see that $Z_\gr(E(\L))$ is just $k$. On the other hand
when the characteristic of $k$ is not two, and if $a=b=1$ and $c$ is
not a root of unity, then we see that $Z_\gr(E(\L))=k[x^2]$ in
$E(\L)$. If $c$ is a primitive $m$-th root of unity, then
$Z_\gr(E(\L))=k[x^2,y^m,z^m]$ in $E(\L)$. By the result of Section
\ref{section:4}, we see that the non-nilpotent elements of
$Z_\gr(E(\L))$ correspond to non-nilpotent elements of
$\HH^*(\L)$. Hence the Hochschild cohomology groups $\HH^n(\L)$ are
non-zero in all even degrees for these choices of $a$, $b$, and $c$.

We leave it to the reader to check that
\[f^n_{u,v,w}=xf^{n-1}_{u-1,v,w}+a^uyf^{n-1}_{u,v-1,t}+
b^uc^vzf^{n-1}_{u,v,w-1}.\]
Furthermore, for $r$ with $0\leq r\leq n$,
\[f^n_{u,v,w}=\sum_{s=\alpha}^{\min\{u,r\}}
\sum_{t=\beta}^{\min\{v,r\}} a^{(u-s)t}b^{(u-s)(r-s-t)}
c^{(v-t)(r-s-t)} f^r_{s,t,r-s-t}f^{n-r}_{u-s,v-t,w+s+t-r},\] where
$\alpha=\max\{0,u+r-n\}$ and $\beta=\max\{0,v+r-n\}$. This formula
yields the elements $c_{pq}(n,i,r)$, and hence gives the structure
constants for the basis of the Koszul dual associated to the elements
$f^n_{u,v,w}$'s and a closed formula for the multiplication in the
Hochschild cohomology ring of $\L$.
\end{example}

\end{document}